\documentclass[10pt]{iopart}
\usepackage{calc}

\usepackage{iopams}
\usepackage{amsthm}
\usepackage[latin1]{inputenc} 
\usepackage{textcomp}

\newcommand{\N}{\nz}
\newcommand{\rz}{\mathbb{R}} 

\newcommand{\grad}{\nabla}

\newcommand{\di}{\,{\rm d}}

\newtheorem{theorem}{Theorem}[section]

\newtheorem{lemma}[theorem]{Lemma}

\makeindex

\usepackage{pstricks}

\newcommand{\text}[1]{{\mathrm{#1}}}
\newcommand{\eqref}[1]{(\ref{#1})}
\newcommand{\bil}{\Gamma}
\newcommand{\R}{\mathbb{R}}
\renewcommand{\N}{\mathbb{N}_0}
\newcommand{\mynegsp}{\vspace{-.5cm}}
\newcommand{\pstpic}{%
\pspicture(0.000000,-0.130000)(0.9470,\pstheight)%
\ifx\nofigs\undefined%
\catcode`@=11%
\psline[linewidth=0.01]{->}(0,0)(0.9470,0)%
\psline[linewidth=0.01]{->}(0,-0.1)(0,\pstheight)%
\rput(-\pstzerox,0){\rput[lt](0.01,0){\rput[lt](\pstonex,-0.01){$1$}}%
\psline[linewidth=0.01](\pstonex,-0.02)(\pstonex,0.02)}%
\rput(0,-\pstzeroy){\rput[rt](-0.01,\pstoney){\rput(-0.01,-0.05){$1$}}%
\psline[linewidth=0.01](-0.02,\pstoney)(0.02,\pstoney)}%
\rput[lt](0,\pstheight){\rput[lt](0,-0.0){%
\begin{minipage}{1cm}
\footnotesize\Mynumb
\end{minipage}}}
}
\newcommand{\pstpend}{%
\catcode`@=12%
\fi%
\endpspicture%
}

\begin{document}
\title[Least Squares Functional for Inverse Sturm-Liouville Problems]{A Least Squares Functional for Solving Inverse Sturm-Liouville Problems}
\author{Norbert Röhrl}
\date{\today}
\address{Fachbereich Mathematik, IADM, University of Stuttgart, 70511 Stuttgart, Germany}
\ead{norbert.roehrl@mathematik.uni-stuttgart.de}
\ams{65L09, 34A55}

\begin{abstract}
  We present a variational algorithm for solving the classical inverse
  Sturm-Liouville problem in one dimension when two spectra are given.
  All critical points of the least squares functional are at global
  minima, which suggests minimization by a (conjugate) gradient
  descent algorithm.  Numerical examples show that the resulting
  algorithm works quite reliable without tuning for particular cases,
  even in the presence of noise. With the right choice of parameters
  the minimization procedure converges very quickly.
\end{abstract}


\section{Introduction}
The one dimensional inverse Sturm-Liouville problem and methods for
its numerical solution are very well studied
\cite{McLaughlin1986,RundellSacks1992,BrownSamkoetal2003}.  Still
there is no general purpose algorithm, which efficiently satisfies all
needs (see eg.~\cite{Schapotschnikow2004}).

We present a variational algorithm which proves to be very robust
under noisy input, does not require any special tuning or additional
input besides the partial spectra, and is quite efficient in the class
of variational methods. In section~\ref{sec:def} we will define the
functional and prove the absence of true local minima, which could
trap our minimization procedure.  Numerical examples, the optimal
choice of the weights and the performance under noisy input will be
discussed in section~\ref{sec:examples}. Section 4 is concerned with
implementation details and compares this algorithm to related
approaches, and section~\ref{sec:proof} finally contains the proof of
the linear independence of the squares of eigenfunctions.

\section{Definition and Properties of the Functional}\label{sec:def}
\noindent We consider the Sturm-Liouville equation
\renewcommand{\theequation}{SL}
\begin{equation}
\label{SL}
{-u'' + q(x)u = \lambda u}
\end{equation}
on $[0,1]$ with $q(x)\in L^2([0,1],\rz)$ real, and separated boundary
conditions
\renewcommand{\theequation}{$\alpha\beta$}
\begin{equation}
 u(0)\cos\alpha +u'(0)\sin\alpha=0,\quad u(1)\cos\beta +u'(1)\sin\beta = 0\,.
\end{equation}
\renewcommand{\theequation}{\arabic{section}.\arabic{equation}}
For the asymptotics 
of the eigenvalues $(\lambda_{n})_{n=0}^\infty$ wrt.~boundary
conditions $(\alpha\beta)$ there are three cases (see eg. \cite{Titchmarsh1962,Isaacsonetal1984,DahlbergTrubowitz1984,PoeschelTrubowitz1987}):
\begin{equation}\label{eigenasymptotics}
\fl \lambda_{n}=\left\{
\begin{array}{ll}
  \pi^2n^2-\frac{2\sin(\beta-\alpha)}{\sin\alpha \sin\beta} + \int_0^1 q \di s + a_n  & \text{ if }\sin\alpha \sin\beta \neq 0\\
  \pi^2(n+1/2)^2+\frac{2\cos\alpha \cos\beta}{\sin(\beta-\alpha)} + \int_0^1 q \di s +a_n &\text{ if } \sin\alpha \sin\beta = 0 \land \sin(\beta-\alpha) \neq 0 \\
  \pi^2(n+1)^2 + \int_0^1 q \di s +a_n &\text{ if } \sin \alpha=\sin \beta=0\\
\end{array}\right.
\end{equation}
where $(a_n)\in l^2$.

It is a classical result, that two spectra determine the
potential uniquely. Fix three angles $\alpha$, $\beta$, and $\gamma$ with
$\sin(\beta-\gamma)\neq0$.  Let $\lambda_{q,i,n}$, $i\in\{1,2\}$ denote the
eigenvalues with respect to the
boundary conditions $(\alpha\beta)$ resp.~$(\alpha\gamma)$.

\begin{theorem}[Borg \cite{Borg1946}, Levinson \cite{Levinson1949}]
\label{BorgLevinson}
Given potentials $Q,q\in L^1[0,1]$ with
\[
\lambda_{Q,i,n}= \lambda_{q,i,n} \qquad \text{ for\ all\ } n\in \N,\quad i=1,2\,,
\] 
then $Q=q$.
\end{theorem}

Suppose we are given (partial) spectral data $\lambda_{Q,i,n}$ with $(i,n)$
in $I \subseteq \{1,2\}× \N$ of an unknown potential $Q$.
For a trial potential $q$ and positive weights $\omega_{i,n}$, we define
the functional
\begin{equation}\label{functional}
G(q) = \sum_{(i,n)\in I} \omega_{i,n}(\lambda_{q,i,n} - \lambda_{Q,i,n})^2 \,.
\end{equation}

If $I$ is infinite and the sequence $(\omega_{i,n})$ is summable, this
converges because of the asymptotics of the eigenvalues \eqref{eigenasymptotics}
\[
G(q) = \sum_{(i,n)\in I} \omega_{i,n}\left(\int_0^1 Q-\int_0^1 q+c_{i,n}\right)^2<\infty \,,
\]
where $c_{i,n}\in l^2, i=1,2$.

If moreover $I=\{1,2\}× \N$, the functional is zero if and only if $q=Q$
(by theorem~\ref{BorgLevinson}). In the case of partial knowledge of
the spectrum, a solution of $G(q)=0$ includes all information given
about the unknown potential $Q$. But we only have a chance to get close to
the original potential if there is not too much information in the
higher eigenvalues, i.e.~if the potential to be recovered is
sufficiently smooth \cite{MarlettaWeikard2005}.

An important question is now, for which sequences $(\lambda_{i,n})$
there is a $q$ with $G(q)=0$, i.e.~$\lambda_{q,i,n} = \lambda_{i,n}$. It is easy
to see with the classical methods \cite{Titchmarsh1962}, that it is 
necessary for the eigenvalues to interlace
\[
   \lambda_{1,n}<\lambda_{2,n}<\lambda_{1,n+1}\qquad \mathrm{ or } \qquad \lambda_{2,n}<\lambda_{1,n}<\lambda_{2,n+1}\,.
\]
It is also known that if we choose either
\begin{enumerate}
\item $\sin \alpha,\sin \beta,\sin \gamma \neq 0$ (Levitan \cite{Levitan1987}) or
\item $\alpha=\beta=0, \sin \gamma\neq 0$ (Dahlberg Trubowitz \cite{DahlbergTrubowitz1984}),  
\end{enumerate}
then the interlacing property in connection with the correct asymptotics
is sufficient for the existence of a potential $q$ with $G(q)=0$. 

An interesting feature of this functional is that all its critical
points are at global minima. To prove this, we first compute the
derivative of $G$. The derivative of $\lambda_{q,i,n}$ wrt.~$q$ in
direction $h$ is
\[
\dot{\lambda}_{q,i,n}[h]= \int_0^1 h g^2_{q,i,n} \di x\,,
\]
where $g_{q,i,n}$ denotes the normalized eigenfunction corresponding to $\lambda_{q,i,n}$ (see \cite{PoeschelTrubowitz1987} for a nice proof).  Thus the
derivative of $G$ is given by
\[
\dot{G}[h](q) = 
2 \sum_{(i,n)\in I} \int_0^1 \omega_{i,n}\left(\lambda_{q,i,n} - \lambda_{Q,i,n} \right)g^2_{q,i,n}h\di x\,. 
\] 
We note that if $n\omega_{i,n}$ is summable, then the gradient
\begin{equation}\label{gradient}
\grad G(q) = 2 \sum_{(i,n)\in I} (n+1)\omega_{i,n}\left(\lambda_{q,i,n} - \lambda_{Q,i,n} \right)\frac1{n+1}g^2_{q,i,n}  
\end{equation}
is in $H^1$ because $\| g^2_{q,i,n}\|_{H^1}=O(n)$ \cite{Titchmarsh1962} and $\lambda_{q,i,n} - \lambda_{Q,i,n}=O(1)$.

Theorem~\ref{linind} below shows that the eigenfunctions $g^2_{q,i,n}$ are linearly
independent in $H^1$. This immediately implies the essential convexity 
of the functional:

\begin{theorem}\label{mainth}
If $I$ is finite or $(n\omega_{i,n})$ is summable, the functional $G$ has no local minima at $q$ with $G(q)>0$, i.e.
\[
\dot{G}[h](q) = 0 \ \forall h\in H^1[0,1] \Longleftrightarrow G(q)=0\,.
\]
\end{theorem}
Thus a conjugate gradient algorithm will not get trapped in local minimas, as we will also observe in the examples.

Finally, we want to address the question, which of the infinitely many
$q$ with $G(q)=0$ our algorithm will select. Let us first look at a
related functional, which actually was also our first try. It is given
by
\begin{equation}\label{altfunc}
\tilde{G}(q) = \sum_{(i,n)\in I} \left(\lambda_{q,i,n} - \lambda_{Q,i,n}+ \int_0^1(Q-q)\di x\right)^2\,, 
\end{equation}
with derivative
\[
\dot{\tilde{G}}[h](q) = 
2 \sum_{(i,n)\in I} \int_0^1 \left(\lambda_{q,i,n} - \lambda_{Q,i,n}
  \int_0^1(Q-q)\di x \right)(g^2_{q,i,n}-1) h\di x \,.
\] 
In case the mean of $Q$ is known, this functional works as well as the
other. Its gradient flow leaves the mean of $q$ constant by
construction. The function $g^2_{q,i,n}-1$ is the gradient of
$\lambda_{q,i,n}-\int_0^1q\di x$ and the direction of strongest increase of
$\lambda_{q,i,n}$ which leaves $\int_0^1q\di x$ fixed.

Choosing $\alpha=\beta=0$ and $\gamma=\pi/2$, the asymptotics of the squared
eigenfunctions are given by
\begin{equation}\label{efasym}
g^2_{q,i,n}=\left\{
\begin{array}{ll}
 1-\cos\big((2n+2)\pi x\big) +O( n^{-1})&, i=1\\
 1-\cos\big((2n+1)\pi x\big) +O( n^{-1})&, i=2\\
\end{array}\right.\,.
\end{equation}
Hence the functions $\{ g^2_{q,i,n}-1| n\in\N, i=1,2\}\cup \{ 1\}$ are almost
orthogonal for large $n$.  
Now, the derivative of our functional \eqref{functional} can be written as
\begin{equation*}
\dot{G}[h](q) = 2 \sum_{(i,n)\in I} \int_0^1 \omega_{i,n}\left(\lambda_{q,i,n} - \lambda_{Q,i,n} \right)\big((g^2_{q,i,n}-1)+1\big)h\di x\,,  
\end{equation*}
and the difference $\lambda_{q+c,i,n}-\int(q+c)\di x=\lambda_{q,i,n}-\int q\di
x$ is invariant under adding a constant function $c$ to $q$.  It follows that
a gradient flow of our functional $G$ leaves $\lambda_{q,i,n}-\int_0^1q\di x$
with $(i,n)\not\in I$ almost invariant for $n$ large enough.

Since also a conjugate gradient descent is just an approximation of
the gradient flow, in practice we do get some little changes in the
higher eigenvalues. A similar argument holds for the case $\sin \alpha,\sin
\beta,\sin \gamma\neq 0$, but there the asymptotics for $i=1,2$ are equal. So we
still have asymptotic orthogonality in $n$, but can not separate $i=1$
and $i=2$.

\section{Numerical Examples}\label{sec:examples}
We use the standard Polak-Ribiere conjugate gradient descent algorithm
\cite{NumericalRecipes} to approximate the gradient flow and thus
minimize the functional. To give the basic idea, we explain the
simpler steepest descent:
\begin{enumerate}
\item choose initial potential $q=q_0$
\item while $G(q_i)$ too big do
  \begin{enumerate}
  \item compute the gradient $\grad G(q_i)$
  \item minimize the one dimensional function $G(q_i - \alpha \grad G(q_i))$ wrt.~$\alpha$  
  \item set $q_{i+1}$ equal to the minimizing potential 
  \end{enumerate}
\end{enumerate}
This straight forward minimization scheme has a major disadvantage:
consecutive gradients are always orthogonal. To avoid this, conjugate
gradient descent computes the direction for the one dimensional
minimization using the current and previous gradients. 

The first natural question is the optimal choice of the constants
$\omega_{i,n}$. A short look at the leading coefficients of the asymptotics
\eqref{eigenasymptotics} may suggest something like
$\omega_{i,n}=(n+1)^{-2}$. But the leading term of $\lambda_{q,i,n} - \lambda_{Q,i,n}$
is $\int_0^1(q-Q) \di x=O(1)$, which corresponds to choosing $\omega_{i,n}=1$
and yields the fastest convergence in all examples.

Figures~\ref{ex1} and~\ref{ex2} were computed with $\omega_{i,n}=1$
resp.~$\omega_{i,n}=(n+1)^{-2}$.  For all plots (except~\ref{ex2} and \ref{ex5}) we chose
the lowest number of iterations after which the plot stabilizes, while
usually already the second iteration reveals the global structure of
the potential. By default we use 30 pairs of eigenvalues, $q_0=0$, $\omega_{i,n}=1$,
and the boundary conditions $\alpha=\beta=0$, $\gamma=\pi/2$.  We give the current
value of the functional $G(q)$ as well as the $L^2$ error
$\Delta_2=\|q-Q\|_2$, and the maximum of the difference of the eigenvalues
$\Delta_\lambda=\max_{(i,n)\in I}\{ |\lambda_{q,i,n} - \lambda_{Q,i,n}|\}$. In addition if there
is a visible difference, we also plot the original potential $Q$ using
a dashed line.
\psset{xunit=8.0cm,yunit=5.0cm}
\newcommand{\pstheight}{0.8}
\newcommand{\pstzeroy}{0.1724}
\newcommand{\pstzerox}{0.1170}
\newcommand{\pstonex}{0.9470}
\newcommand{\pstoney}{0.3492}
\newcommand{\numberofit}{\text{iter\#}}
\newcommand{\Mynumb}{%
\[\fl\arraycolsep0.1cm\begin{array}{rcl}
    G(q) &=& 4.1589\cdot10^{-5}\\
    \Delta_2 &=& 0.182951\\
    \Delta_\lambda &=& 2.2881\cdot10^{-3}\\
    \numberofit &=& 5
  \end{array}\]}
\input{gnuplothead}
\begin{figure}
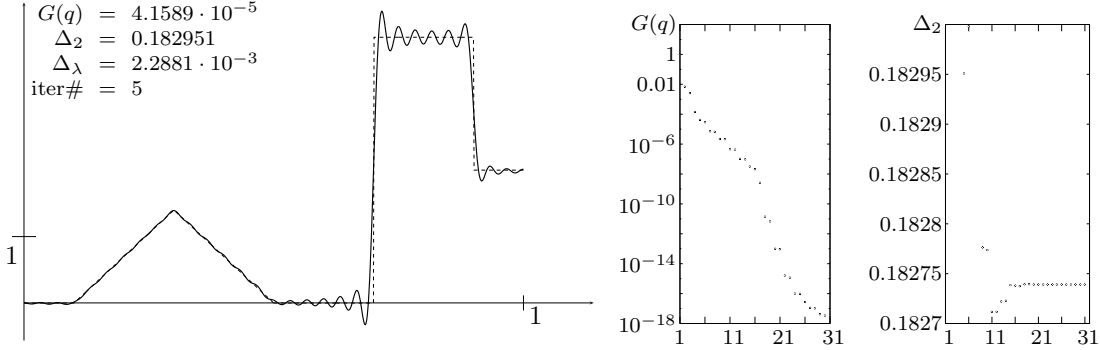

\begin{tabular}{ccc}
%
\pstpic\rput(-\pstzerox,-\pstzeroy){{\input{val2_n30}}{\input{val2_func}}}\pstpend& 
\phantom{---}\psset{xunit=2.5cm,yunit=4.5cm}\footnotesize \input{convG_n30_e0} &
\phantom{-----}\psset{xunit=2.5cm,yunit=4.5cm}\footnotesize \input{convl2_n30_e0} \\
\end{tabular}\mynegsp
\caption{%
Plot using $\omega_{i,n}=1$. The right two graphs show
$G(q)$ \\and $\Delta_2$ over the number of iterations.
}
\label{ex1}
\end{figure}
\begin{figure}
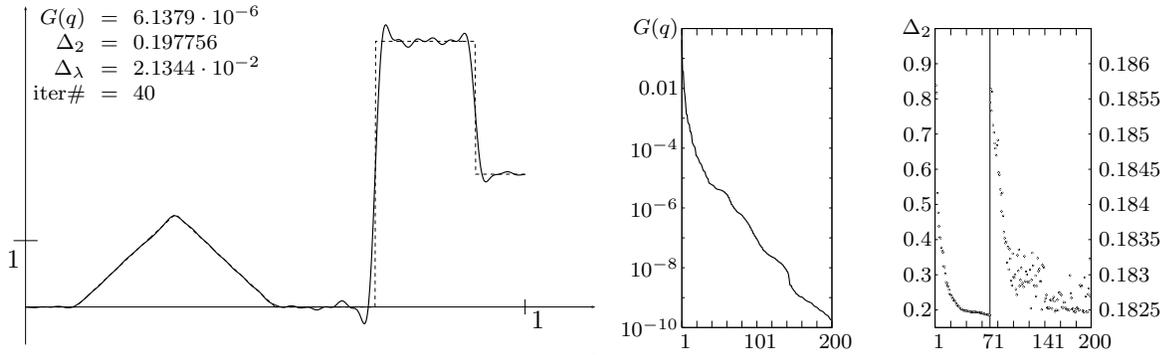

\begin{tabular}{ccc}
\renewcommand{\Mynumb}{%
\[\fl\arraycolsep0.1cm\begin{array}{rcl} G(q) &=& 6.1379\cdot10^{-6}\\\Delta_2 &=&0.197756 \\
\Delta_\lambda &=& 2.1344\cdot10^{-2}\\\numberofit &=& 40\end{array}\]}
%
\pstpic\rput(-\pstzerox,-\pstzeroy){{\input{val2_n30_e2}}{\input{val2_func}}}\pstpend
& 
\phantom{---}\psset{xunit=2.5cm,yunit=4.5cm}\footnotesize \input{convG_n30_e2} &
\phantom{-----}\psset{xunit=2.5cm,yunit=4.5cm}\footnotesize \input{convl2_n30_e2} \\
\end{tabular}\mynegsp
\caption{Plot using $\omega_{i,n}=(n+1)^{-2}$.  On the right side of the  $\Delta_2$ plot we zoom in to show the convergence process in more detail.}
\label{ex2}
\end{figure}
We see that both $G(q)$ and $\Delta_2$ converge much faster for
$\omega_{i,n}=1$. With $\omega_{i,n}=(n+1)^{-2}$ on the other hand, $q$ moves
through smoother functions (as in figure~\ref{ex2}). In contrast to
the optical impression these are in general worse approximations
wrt.~$G(q)$ and $\Delta_2$. But finally $q$, after around 100
iterations, will also converge to the function shown in figure
\ref{ex1}. From the good average convergence of $G(q)$ in both cases,
we can tell, that our algorithm also in practice does not get near a
local minimum, as proved in theorem~\ref{mainth}.
  
The optimal choice of boundary conditions can already be guessed from the
asymptotics of the squared eigenfunctions \eqref{efasym}.  Choosing
$\alpha=\beta=0$ and $\gamma=\pi/2$, the functions $g_{q,i,n}^2-1$ are almost
orthogonal for large $n$ and $i=1,2$ (c.f.~\eqref{efasym}). In contrast, if
$\sin\alpha$, $\sin\beta$, and $\sin\gamma$ are all non zero, $g_{q,1,n}^2-1$
and $g_{q,2,n}^2-1$ will get close for large $n$.  Therefore, in the latter
case, the algorithm will converge much slower. This can also clearly be
observed in numerical examples, like in figure~\ref{ex5}. But also there,
after approximately 620 iterations, we finally will get a solution which is
close to figure~\ref{ex1}.

\begin{figure}
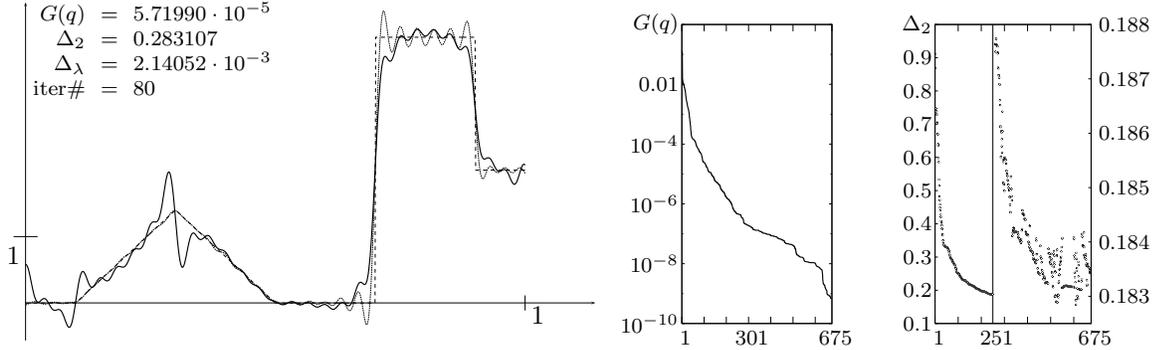

\begin{tabular}{ccc}
\renewcommand{\Mynumb}{%
\[\fl\arraycolsep0.1cm\begin{array}{rcl} G(q) &=& 5.71990\cdot 10^{-5}\\\Delta_2 &=&0.283107\\
\Delta_\lambda &=& 2.14052\cdot10^{-3}\\\numberofit &=& 80\end{array}\]}
\pstpic\rput(-\pstzerox,-\pstzeroy){{\input{val2_bd2_n30}}{\input{val2_func}}}\pstpend
& 
\phantom{---}\psset{xunit=2.5cm,yunit=4.5cm}\footnotesize \input{invbdconvG} &
\phantom{-----}\psset{xunit=2.5cm,yunit=4.5cm}\footnotesize \input{convl2bd1} \\
\end{tabular}\mynegsp
\caption{Plot with bounday conditions $\alpha=\beta=\pi/4$, $\gamma=-\pi/4$.\\ The light
  graph shows the 620th iteration, where the iteration\\ stabilizes.}
\label{ex5}
\end{figure}

\begin{figure}
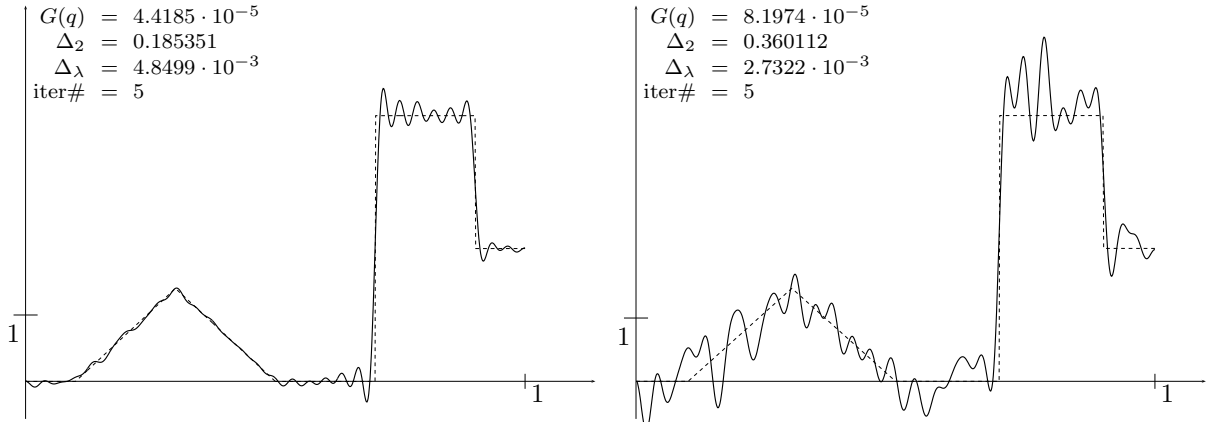

\begin{tabular}{cc}
\renewcommand{\Mynumb}{%
\[\fl\arraycolsep0.1cm\begin{array}{rcl} G(q) &=& 4.4185\cdot10^{-5}\\\Delta_2 &=& 0.185351\\
\Delta_\lambda &=& 4.8499\cdot10^{-3}\\\numberofit &=& 5\end{array}\]}

%
\pstpic\rput(-\pstzerox,-\pstzeroy){{\input{val2_n30_r0.01}}{\input{val2_func}}}\pstpend &
\renewcommand{\pstzerox}{0.0850}%
\renewcommand{\pstzeroy}{0}%
\renewcommand{\pstoney}{0.1684}%
\renewcommand{\pstonex}{0.9470}%
\renewcommand{\Mynumb}{%
\[\fl\arraycolsep0.1cm\begin{array}{rcl} G(q) &=& 8.1974\cdot10^{-5}\\\Delta_2 &=& 0.360112\\
\Delta_\lambda &=& 2.7322\cdot10^{-3}\\\numberofit &=& 5\end{array}\]}

%
\pstpic\rput(-\pstzerox,-\pstzeroy){{\input{val2_n30_r0.1}}{\input{val2_func6}}}\pstpend\\
\end{tabular}\mynegsp
\fl\caption{Plots with added white noise
of absolute value smaller \\ or equal $0.01$ resp.~$0.1$.}
\label{ex3}
\end{figure}

Another important aspect for applications is stability against noise
in the given spectral data. In figure~\ref{ex3} we computed two
examples with random noise $|\tilde{\lambda}_{Q,i,n} - \lambda_{Q,i,n}|\leq r$, with
$r=0.01$ and $r=0.1$, respectively. It is remarkable, that the
convergence speed in $G(q)$ is not significantly affected by the
random noise; both examples reach $G(q)\sim10^{-18}$ in 30 iterations.
\renewcommand{\pstheight}{1}

For testing the robustness, we also fed the functional with the following
random sequence
\[\fl
9.99742, 11.6265, 14.4527, 23.9247, 26.2413, 31.091, 40.6658, 48.1088, 53.5093,
60.9088, 
\]
which is far from the generic asymptotics. Convergence of $G(q)$ is
slow but steady. In log scale the graph of $G(q)$ over the number of
iteration looks similar to those given above. From $G(q)=36287$
initially, we get down to $G(q)=3.56528\cdot 10^{-9}$ in 450 iterations.

\begin{figure}
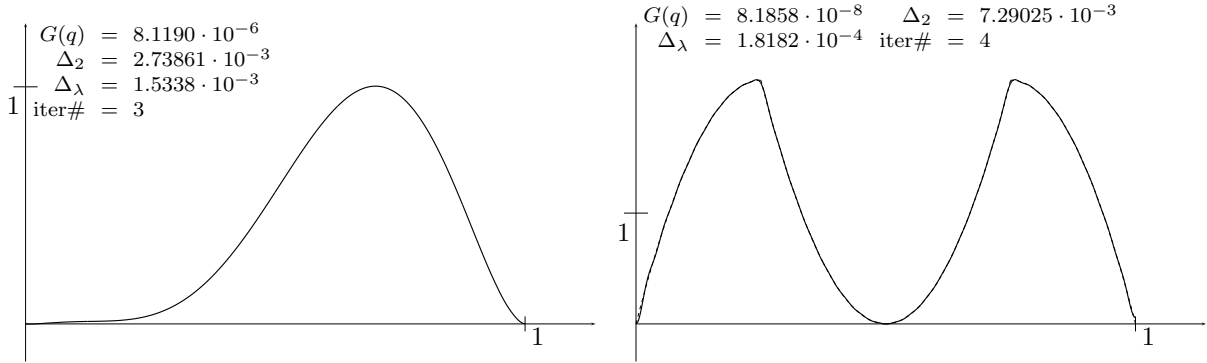

\renewcommand{\pstheight}{0.8}
\begin{tabular}{cc}
\renewcommand{\Mynumb}{%
\[\fl\arraycolsep0.1cm\begin{array}{rcl} G(q) &=& 8.1190\cdot10^{-6}\\\Delta_2 &=& 2.73861\cdot10^{-3}\\
\Delta_\lambda &=& 1.5338\cdot10^{-3}\\\numberofit &=& 3\end{array}\]}
\renewcommand{\pstoney}{0.8417}%
\renewcommand{\pstzerox}{0.1170}%
\renewcommand{\pstzeroy}{0.2103}%
%
\pstpic\rput(-\pstzerox,-\pstzeroy){{\input{val1_n30_e0}}{\input{empty}}}\pstpend&
\renewcommand{\Mynumb}{%
\rput[lt](0,0.25cm){\begin{minipage}{2cm}{\[\arraycolsep0.1cm\fl\begin{array}{rclrcl} G(q) &=& 8.1858\cdot10^{-8}&\Delta_2 &=& 7.29025\cdot10^{-3}\\
\Delta_\lambda &=& 1.8182\cdot10^{-4}&\numberofit &=& 4\end{array}\]}
\end{minipage}}}
\renewcommand{\pstoney}{0.5260}%
\renewcommand{\pstzerox}{0.1170}%
\renewcommand{\pstzeroy}{0.2313}%
\pstpic\rput(-\pstzerox,-\pstzeroy){{\input{val3_n30_e0}}{\input{val3_testfun}}}\pstpend \\
\end{tabular}\mynegsp
\caption{Continuous examples.}
\label{ex4}
\end{figure}

Finally figure~\ref{ex4} shows the results for the other examples of
\cite{BrownSamkoetal2003,RundellSacks1992}.
\section{Implementation and Comparison}
The implementation of the conjugate gradient minimization algorithm is
straight forward (see eg.~\cite{NumericalRecipes}), but our
application crucially depends on a reliable eigenvalue problem solver.
After some trial and error we settled for d02kdf from the Numerical
Algorithms Group (NAG). For standard numerical routines we used the
GNU Scientific Library.

Another important point is the internal representation of the
functions.  We used cubic splines on 2000 intervals of equal length.
For 30 pairs of eigenvalues we already get satisfactory results at 180
intervals. But since the computation time depends only roughly linear
on this number, we did not study its influence on the quality of the
approximation.

Because all high level routines, like the eigenvalue solver, only sample
these functions on some points, the details of the representation are
not that important. But since very many values are needed, this is the
part in the algorithm, where speed optimizations are most valuable.

Compared to the variational algorithm by Brown et
al.~\cite{BrownSamkoetal2003}, which uses other spectral data, our
functional is more expensive to compute, because they only have to
solve initial value problems, whereas the present algorithm requires
solution of the eigenvalue problems. On the other hand they often need
about 10-100 times as many iterations to get similar results. Another
plus on our side is, that even in presence of noise, we know to get
closer to our goal each step and do not have to regularize the process
by limiting the number of iterations.

The method of Rundell and Sacks \cite{RundellSacks1992} uses the same
spectral data as our algorithm and is without any doubt much faster,
but needs the mean $\int_0^1Q \di x$ as additional input, which has to be
guessed from the (partial) spectral data.  Our method does this
automatically. In the case of figure~\ref{ex1} for example, the error
of the mean is $2\cdot10^{-3}\%$ and using 5 pairs of eigenvalues we still
only get an error of about $0.2\%$. This is probably better than an
independent algorithm for computing the mean value from spectral data
could do.

\section{Linear Independence of Squared Eigenfunctions }\label{sec:proof}
Define the Wronskian
 $[f,g]= fg'-f'g$ and the bilinear form
\[
\begin{array}{rrcl}
\bil: &H^1([0,1],\R)^2& \longrightarrow & \R\\
   &   (f,g) &\mapsto &  \int_0^1 [f,g] \di x
\end{array}
\]
which is bounded by 
\[
|\bil(f,g)| \leq \|f\|_{H_1} \|g\|_{H_1},\text{ i.e. } \|\bil(f,\cdot) \| = \|f\|_{H_1}\,. 
\]
(We use the definition $\|f\|_{H^1}=\sqrt{\|f\|_{L^2}^2+ \|f'\|_{L^2}^2}$ with
distributional derivatives.) 
In particular
$\bil$ is continuous on $H^1$.  We have the following rules for the Wronskian:
\begin{enumerate}
\item $[fg,hj] = gj[f,h] + fh[g,j] = fj[g,h] + gh[f,j]$
\item If the functions $f_1$ and $f_2$ fulfill the condition 
\[
f_i(a)\cos\alpha +f_i'(a)\sin\alpha=0 \qquad i = 1,2\,,
\]
then $[f_1,f_2](a) = 0$.
\item
For two arbitrary solutions $f_1$ and $f_2$ of the equation
\eqref{SL} with eigenvalue parameters $\lambda_1$ and $\lambda_2$ we
have
\[
[f_1,f_2]' = f_1f_2'' - f_2f_1'' = ( \lambda_1-\lambda_2 ) f_1f_2\,.
\]
\end{enumerate}



\noindent
In this section we are only talking about a single $q$ and therefore
drop it from the index. Let $s_{i,n}$ and $c_{i,n}$ be the solutions
of the differential equation \eqref{SL} for the eigenvalue parameter
$\lambda_{i,n}$ and initial values
\[
\begin{array}{rclcrcl}
  s_{i,n}(1) &=& \sin \beta\,,&\phantom{---}&\qquad c_{i,n}(1) &=& \sin \gamma\,, \\
  s_{i,n}'(1) &=& -\cos \beta\,,& &\qquad c_{i,n}'(1) &=& -\cos \gamma\,.
\end{array}
\]
It is well known and easy to see that the Wronskian of these solutions 
is constant and we can compute its value at 1:
\[
[s_{i,n},c_{i,n}] = 
-\sin \beta \cos \gamma + \cos\beta \sin \gamma = 
\sin(\gamma -\beta). 
\] 
The normalized eigenfunctions are $g_{1,n}=s_{1,n}/\|s_{1,n}\|_2$ and
$g_{2,n}=c_{2,n}/\|c_{2,n}\|_2$. Now we prove the central lemma, which
builds on the ideas of a similar result for the Dirichlet case in the
book of Pöschel and Trubowitz\cite{PoeschelTrubowitz1987}. A result in
the same spirit can also already be found in the paper of Borg
\cite{Borg1946}.
\begin{lemma}
Given three angles $\alpha$, $\beta$, and $\gamma$ with $\sin(\beta-\gamma)\neq0$ and denote the
$L^2$-normalized eigenfunctions of the Sturm-Liouville equation
with boundary conditions $(\alpha\beta)$ and $(\alpha\gamma)$ by $g_{i,n}$, $i=1,2$. With
$s_{i,n}$, $c_{i,n}$ as defined above we have the following relations
for the squared eigenfunctions for all $i,j\in \{1,2\}$ and $m,n\in\N$: 
\begin{enumerate}
  \item $\displaystyle \bil(g_{i,n}^2,g_{i,m}^2) = 0$ 
\item $\displaystyle \bil(c_{i,n} s_{i,n},g_{j,m}^2) = (-1)^i \sin(\gamma -\beta)\delta_{n,m}\delta_{i,j}$ 
\end{enumerate}
\end{lemma}
\begin{proof}
i)
\[
\bil(g_{i,n}^2,g_{i,m}^2) = 2 \int_0^1 g_{i,n}g_{i,m}[g_{i,n},g_{i,m}] \di x
\]
If $n=m$, the term clearly vanishes. 
If $n\neq m$, we get 
\[
\fl
\bil(g_{i,n}^2,g_{i,m}^2) = \frac2{ \lambda_{i,n}-\lambda_{i,m} } 
\int_0^1 [g_{i,n},g_{i,m}]'[g_{i,n},g_{i,m}] \di x = \frac1{ \lambda_{i,n}-\lambda_{i,m} }[g_{i,n},g_{i,m}]^2\big|_0^1 =0\,,
\]
because $g_{i,n}$ and $g_{i,m}$ satisfy the same boundary conditions (rule ii). 

\noindent
ii)
\[
\bil(c_{i,n}s_{i,n},g_{j,m}^2) = \int_0^1 \big(c_{i,n}g_{j,m}[s_{i,n},g_{j,m}] + s_{i,n}g_{j,m}[c_{i,n},g_{j,m}]\big) \di x
\]
If $i=j=1, m=n$, the first term vanishes and we are left with
\[
\int_0^1 s_{1,m}g_{1,m}[c_{1,m},g_{1,m}] = \int_0^1 g^2_{1,m}[c_{1,m},s_{1,m}] = -\sin(\gamma -\beta)\,,
\]
and if $i=j=2, m=n$, the second term vanishes and we get
\[
\int_0^1 c_{2,m}g_{2,m}[s_{2,m},g_{2,m}] = \int_0^1 g^2_{2,m}[s_{2,m},c_{2,m}] = +\sin(\gamma -\beta)\,.
\]
If $(i,n)\neq(j,m)$ we compute 
\begin{eqnarray*}\fl
\bil(c_{i,n}s_{i,n},g_{j,m}^2) =& \frac1{ \lambda_{i,n}-\lambda_{j,m} }\int_0^1 \big([c_{i,n},g_{j,m}]'[s_{i,n},g_{j,m}] + [s_{i,n},g_{j,m}]'[c_{i,n},g_{j,m}]\big) \di x = \\
 & \frac1{ \lambda_{i,n}-\lambda_{j,m} } [s_{i,n},g_{j,m}][c_{i,n},g_{j,m}] \big|_0^1 =0 
\end{eqnarray*}
We note that $\lambda_{i,n}-\lambda_{j,m}\neq 0$ in this case, because the eigenfunctions satisfy the same boundary conditions at $0$ and different boundary conditions at $1$.
\end{proof}

\begin{theorem}\label{linind}
With the notations of the above lemma, the set of squared eigenfunctions
\[
 \big\{ g_{i,n}^2 | (i,n)\in \{1,2\} × \N\big\} 
\]
is linearly independent in $H^1$.
\end{theorem}
\begin{proof}
Suppose for some fixed $(i,n)$ we have
\[
g^2_{i,n}=\sum_{k\in \N} a_kg^2_k
\]
in $H^1$,
where $a_k\in\rz$ and $g_k=g_{j_k,m_k}$ with $(j_k,m_k)\neq(i,n)$. But
this would imply
\[\fl
(-1)^i\sin(\gamma -\beta) = \bil (c_{i,n}s_{i,n},g_{i,n}^2) = 
\bil \left(c_{i,n}s_{i,n},\sum_{k\in\N} a_kg^2_k\right) = \sum_{k\in\N}\bil (c_{i,n}s_{i,n},a_kg^2_k) =0 
\] 
\end{proof}

\ack
The author wishes to thank Ian
Knowles and Rudi Weikard for drawing his attention to this topic and
for the interesting discussions during his time at the University of
Alabama at Birmingham.

\section*{References}

\bibliographystyle{abbrv}
\bibliography{inp}

\end{document}